\documentclass{amsart}
\usepackage{amsmath, amscd, amssymb, amsthm}
\usepackage{bbm}
\usepackage{latexsym}
\usepackage{amsfonts}
\usepackage{graphicx}
\usepackage[all,cmtip]{xy}

\usepackage[
colorlinks=true, linkcolor=black,breaklinks=true,
urlcolor=blue, citecolor=black]{hyperref}%
\usepackage{geometry}
\newtheorem{theorem}{Theorem}
\newtheorem{lemma}{Lemma}

\newtheorem{question}[theorem]{Question}

\newtheorem{definition}{Definition}

\newcommand{\be}{\begin{equation}}
\newcommand{\ee}{\end{equation}}
\newcommand{\bea}{\begin{eqnarray}}
\newcommand{\eea}{\end{eqnarray}}

\def\XXint#1#2#3{{\setbox0=\hbox{$#1{#2#3}{\int}$ }
\vcenter{\hbox{$#2#3$ }}\kern-.6\wd0}}

\begin{document}

\title[Flat Hermitian Lie algebras are K\"ahler]{Flat Hermitian Lie algebras are K\"ahler}

\author{Dongmei Zhang}
\address{Dongmei Zhang. School of Mathematical Sciences, Chongqing Normal University, Chongqing 401331, China}
\email{{2250825921@qq.com}}\thanks{Zhang is supported by Chongqing graduate student research grant No.\,CYB240227. The corresponding author Zheng is partially supported by National Natural Science Foundations of China with the grant No.\,12141101 and 12471039,  Chongqing Normal University grant 24XLB026, and is supported by the 111 Project D21024.}

\author{Fangyang Zheng}
\address{Fangyang Zheng. School of Mathematical Sciences, Chongqing Normal University, Chongqing 401331, China}
\email{20190045@cqnu.edu.cn; franciszheng@yahoo.com} \thanks{}

\subjclass[2020]{53C55 (primary), 53C05 (secondary)}
\keywords{flat Hermitian manifolds; flat Lie groups; Hermitian Lie algebras; Levi-Civita flat; K\"ahler flat}

\begin{abstract}
In 1976, Milnor classified all Lie groups admitting a flat left-invariant metric. They form a special type of unimodular 2-step solvable groups. Considering Lie groups with  Hermitian structure, namely, a left-invariant complex structure and a compatible left-invariant metric, in 2006, Barberis-Dotti-Fino obtained among other things full classification of all Lie groups with Hermitian structure that are K\"ahler and flat. In this note, we examine Lie groups with a Hermitian structure that are flat, and show that they actually must be K\"ahler, or equivalently speaking, a flat Hermitian Lie algebra is always K\"ahler. In the proofs we utilized analysis on the Hermitian geometry of 2-step solvable Lie groups developed by Freibert-Swann and by Chen and the second named author.
\end{abstract}

\maketitle


\markleft{Zhang and Zheng}
\markright{Flat Hermitian Lie algebras are K\"ahler}

\section{Introduction and statement of result}

In 1976, Milnor  studied the curvature properties of Lie groups with left-invariant metrics. His \cite[Theorem 1.5]{Milnor} states that, {\em a Lie group with left-invariant  metric is flat if and only  if the associated Lie algebra splits as an orthogonal direct sum ${\mathfrak h} \oplus {\mathfrak u}$ where ${\mathfrak h}$ is a commutative subalgebra, ${\mathfrak u}$ is a commutative ideal, and where  the linear transformation $\mbox{ad}_x$ is skew-adjoint for every $x\in {\mathfrak h}$. } In particular, such Lie groups are always unimodular and 2-step solvable (or abelian).

Let $G$ be a Lie group and ${\mathfrak g}$ its Lie algebra. As is well-known, left-invariant metrics on $G$ are in one-one correspondence with metrics (i.e., inner products) on ${\mathfrak g}$. Similarly, left-invariant complex structures on $G$ are in one-one correspondence with complex structures on ${\mathfrak g}$, namely, linear maps $J: {\mathfrak g} \rightarrow {\mathfrak g}$ satisfying $J^2=-I$ and the integrability condition:
\begin{equation} \label{integrability}
 [x,y] -[Jx,Jy] + J[Jx,y] + J[x,Jy] =0, \ \ \ \ \ \ \forall \ x,y\in {\mathfrak g}.
 \end{equation}

In 2006, Barberis, Dotti and Fino strengthened Milnor's theorem and gave the following more detailed description of flat Lie groups. See  \cite[Props 2.1 and  2.2]{BDF} for more details, here we rephrased it slightly in the style of \cite[Appendix]{VYZ}.
\begin{theorem}[\cite{BDF}]
Let $({\mathfrak g}, g)$ be a Lie algebra with a metric. Then $g$ is flat if and only if ${\mathfrak g}$ is the orthogonal direct sum ${\mathfrak h} \oplus {\mathfrak z} \oplus {\mathfrak g}'$, where ${\mathfrak z}$ is the center of ${\mathfrak g}$, ${\mathfrak g}'=[{\mathfrak g}, {\mathfrak g}]$ is the commutator, ${\mathfrak h}$ is an abelian subalgebra, satisfying: 
\begin{enumerate}
\item $\dim_{\tiny \mathbb R}({\mathfrak g}')=2p$ is even, with an orthonormal basis $\{ \varepsilon_1, \ldots , \varepsilon_{2p}\}$  so that
\item $[x,\varepsilon_{2i-1}]=f_i(x)\varepsilon_{2i}$, $[x,\varepsilon_{2i}]=-f_i(x)\varepsilon_{2i-1}$, $\forall\ 1\leq i\leq p$ and $\forall \ x\in {\mathfrak h}$, where $f=(f_1, \ldots , f_p)$ is an injective linear map from ${\mathfrak h}$ into ${\mathbb R}^{p}$ with each $f_i\not\equiv 0$. 
\end{enumerate}
\end{theorem}
Note that in \cite[Appendix]{VYZ} the condition `each $f_i\not\equiv 0$' was missing, which is needed as when $f_i\equiv 0$, $\varepsilon_{2i-1}$ and $\varepsilon_{2i}$  are not  in ${\mathfrak g}'$.

Barberis, Dotti and Fino showed that, when $\dim_{\tiny \mathbb R}({\mathfrak g})$ is even, a flat Lie algebra $({\mathfrak g}, g)$ always admit a compatible complex structure $J$ so that $g$ is K\"ahler. Conversely, their \cite[Cor 2.1, Prop 3.1]{BDF} gives a complete classification of K\"ahler flat Lie algebras (see also \cite[Theorem A.2]{VYZ} in the Appendix for the slight reformulation here):

\begin{theorem}[\cite{BDF}] \label{thm2}
Let $({\mathfrak g}, J, g)$ be a Lie algebra with a Hermitian structure. Then $g$ is K\"ahler flat if and only if ${\mathfrak g}$ is the orthogonal direct sum ${\mathfrak h} \oplus {\mathfrak z} \oplus {\mathfrak g}'$, where ${\mathfrak z}$ is the center of ${\mathfrak g}$, ${\mathfrak g}'=[{\mathfrak g}, {\mathfrak g}]$ is the commutator, ${\mathfrak h}$ is an abelian subalgebra, satisfying: 
\begin{enumerate}
\item $J{\mathfrak g}'={\mathfrak g}'$, with an orthonormal basis $\{ \varepsilon_1, \ldots , \varepsilon_{2p}\}$  so that $\varepsilon_{2i}=J\varepsilon_{2i-1}$ $\forall\ 1\leq i\leq p$, and
\item $[x,\varepsilon_{2i-1}]=f_i(x)\varepsilon_{2i}$, $[x,\varepsilon_{2i}]=-f_i(x)\varepsilon_{2i-1}$, $\forall\ 1\leq i\leq p$ and $\forall \ x\in {\mathfrak h}$, where $f=(f_1, \ldots , f_p)$ is an injective linear map from ${\mathfrak h}$ into ${\mathbb R}^{p}$ with each $f_i\not\equiv 0$, and
\item there are further orthogonal decomposition ${\mathfrak h} = ({\mathfrak h}\cap J{\mathfrak h})\oplus {\mathfrak h}_1$ and ${\mathfrak z} = ({\mathfrak z}\cap J{\mathfrak z})\oplus {\mathfrak z}_1$, so that $J$ restricts to an isomorphism between ${\mathfrak h}_1$ and ${\mathfrak z}_1$.  
\end{enumerate}
\end{theorem}

Now suppose $({\mathfrak g}, J, g)$ is a Lie algebra with a Hermitian structure such that $g$ is flat. It is natural to wonder about the classification of all flat Hermitian Lie algebras, and in particular, if such a $g$ must be K\"ahler or not? The main purpose of this short article is just to answer this question, and it turns out that such a $g$ is always K\"ahler, whose classification was given by Theorem \ref{thm2}. 

\begin{theorem} \label{thm3}
Let $({\mathfrak g}, J, g)$ be a Lie algebra with a Hermitian structure such that $g$ is flat. Then $g$ is K\"ahler, so the structure is given by Theorem \ref{thm2}.
\end{theorem}



\noindent {\bf Remark.}  (1). First of all, there are plenty of flat Hermitian manifolds that are not K\"ahler, even in the compact case. For instance, for each $n\geq 3$, there are examples of compact Hermitian $n$-manifolds that are flat but non-K\"ahler. In other words, on some flat tori $(T^{2n}_{\mathbb R},g_0)$ there exist (plenty of) complex structures $J$ compatible with $g_0$ such that $g_0$ is not K\"ahler with respect to $J$ (hence $(T^{2n}_{\mathbb R},J)$ is not a complex torus). For $n=3$, all such $J$ were classified by Khan, Yang, and the second named author in \cite{KYZ}, but for $n\geq 4$, such a classification is still unknown. For $n=2$, all compact flat Hermitian surfaces are K\"ahler, although there are examples of  complete, flat Hermitian surfaces which are not K\"ahler. The above Theorem \ref{thm3} says that for Lie groups (with left-invariant complex structures and metrics), flat Hermitian metrics are always K\"ahler. 

\vspace{0.1cm}

(2). Given any flat Lie algebra $({\mathfrak g},g)$ with $\dim_{\tiny \mathbb R}({\mathfrak g})=2n$ being even,  there are plenty of  almost complex structures $J$ on ${\mathfrak g}$ compatible with the metric $g$. In fact, any skew-symmetric, orthogonal $2n\times 2n$ matrix will provide us with such a $J$. Of course such a $J$ does not need to preserve ${\mathfrak g}'$ in general. In order for $J$ to be a complex structure, the only condition it needs to satisfy is the integrability condition (\ref{integrability}). Just from the algebraic point of view, it seems unlikely that this single condition could rule out all those $J$ with $J{\mathfrak g}' \neq {\mathfrak g}'$. Yet Theorem \ref{thm3} says that is indeed the case, which illustrates the beauty of geometry once again.

\vspace{0.1cm}

(3). For Hermitian Lie algebras, being K\"ahler is a strong restriction. For instance, the classic result of Hano \cite{Hano} states that any unimodular K\"ahler Lie algebra is always flat.

\vspace{0.15cm}

In the past a couple of decades, the Hermitian geometry on Lie groups (or equivalently, Hermitian Lie algebras) has been  studied by many authors, such as Gray, Salamon, Ugarte, Fino, Vezzoni, Tomassini, Podesta, Grantcharov, Kasuya, Andrada, Barberis, Dotti, Angella,  Di Scala, Lauret, Lafuente, Arroyo, Nicolini, and others. As a small sampler, see \cite{ABD, AOUV, AL, AN, CFGU, DLV, EFV, FG, FK, FKV, FinoT,  GiustiPodesta, LZ} and the references therein.

The proof of Theorem \ref{thm3} is based on Hermitian geometric properties for 2-step solvable Lie algebras. In \cite{FS, FS2}, Freibert and Swann give a systematic treatment on the Hermitian geometry of 2-step solvable Lie algebras. In particular, they confirmed the Fino-Vezzoni Conjecture  for all such Hermitian Lie algebras of 
the pure types, which means either $J{\mathfrak g}'\cap {\mathfrak g}'=0$, or $J{\mathfrak g}'= {\mathfrak g}'$, or $J{\mathfrak g}'+ {\mathfrak g}'={\mathfrak g}$. Fino-Vezzoni Conjecture (\cite{FV, FV2}) states that if a compact complex manifold admits a balanced metric and a pluriclosed metric, then it must admit a K\"ahler metric. 

In a recent work \cite{CZ}, Chen and the second named author confirmed the Streets-Tian Conjecture for all 2-step solvable Lie algebras. The conjecture (\cite{ST, ST2}) states that if a compact complex manifold admits a Hermitian symplectic metric, then it must admit a K\"ahler metric. Results from \cite{FS, FS2} and \cite{CZ} on the Hermitian structures of 2-step solvable Lie algebras are utilized in the proof of Theorem \ref{thm3}.

\vspace{0.3cm}

\section{Preliminaries }

Let $G$ be a connected, simply-connected, even-dimensional Lie group, and ${\mathfrak g}$ its Lie algebra. Let $(J,g)$ be a left-invariant Hermitian structure on $G$, and we denote by the same letters $(J,g)$ for the corresponding Hermitian structure on ${\mathfrak g}$. For convenience, we will write $g=\langle , \rangle$ and extend it bi-linearly over ${\mathbb C}$. By a {\em frame}  on ${\mathfrak g}$ we mean a basis $\{ e_1, \ldots , e_n\}$ of the complex vector space ${\mathfrak g}^{1,0}=\{ x-\sqrt{-1}Jx \mid x\in {\mathfrak g}\}$. Similarly, a {\em coframe}  on ${\mathfrak g}$ means a basis $\{ \varphi_1 , \ldots , \varphi_n\}$ of the dual vector space  $({\mathfrak g}^{1,0})^{\ast}$. The coframe $\varphi$ is said to be dual to $e$ if $\varphi_i(e_j)=\delta_{ij}$ for any $1\leq i,j\leq n$. 

Let us fix a frame $e$ on ${\mathfrak g}$ and let $\varphi$ be its dual cofarme. Note that here we do not assume that $e$ to be unitary. Following the notations of \cite{VYZ,YZ} and \cite{CZ}, let us denote by
\begin{equation}  \label{CD1}
C_{ik}^j = \varphi_j( [e_i,e_k] ) , \ \ \ \ \ D_{ik}^j = \overline{\varphi}_i ( [\overline{e}_j, e_k]), \ \ \ \ \ 1\leq i,j,k\leq n,
\end{equation}
for the structure constants. This is equivalent to the following
\begin{equation}
[e_i,e_j] = \sum_{k=1}^n C_{ij}^k e_k , \ \ \ \ \ [e_i,\overline{e}_j ] = \sum_{k=1}^n \big( \overline{D^i_{kj}} \,e_k -  D_{ki}^j \overline{e}_k \big), \ \ \ \ \ 1\leq i,j\leq n.   \label{CD2}
\end{equation}
We can extend each $e_i$  to a left-invariant vector field on $G$ and still denote it by $e_i$. Thus $e$ becomes a global  frame of type $(1,0)$ tangent vector fields on $G$ as a Hermitian manifold. Under this frame, the first structure equation and the first Bianchi identity become
\begin{equation}
d\varphi_{i} = -\frac{1}{2} \sum_{j,k=1}^n  C^i_{jk}\,\varphi_j \wedge \varphi_k - \sum_{j,k=1}^n \overline{D^j_{ik}} \,\varphi_{j} \wedge \overline{\varphi}_{k}, \ \ \ \ \ \ \forall \ 1\leq i\leq n.  \label{SED}
\end{equation}
\begin{equation}
\left\{   \begin{split}
 \ \sum_{t=1}^n \big( C^t_{ij}C^{\ell}_{tk} + C^t_{jk}C^{\ell}_{ti} + C^t_{ki}C^{\ell}_{tj} \big) \ = \ 0, \hspace{3.13cm} \\
 \ \sum_{t=1}^n \big( C^t_{ik}D^{\ell}_{jt} + D^t_{ji}D^{\ell}_{tk} - D^t_{jk}D^{\ell}_{ti} \big) \ = \ 0 , \hspace{2.92cm}   \\
 \ \sum_{t=1}^n \big( C^t_{ik}\overline{D^{t}_{j\ell }} - C^j_{tk}\overline{D^{i}_{t\ell }} + C^j_{ti}\overline{D^{k}_{t\ell }}  - D^{\ell}_{ti}\overline{D^{k}_{jt }} + D^{\ell}_{tk} \overline{ D^{i}_{jt }}  \big) \ = \ 0,  
\end{split}  \right. \label{CCCD}
\end{equation}
for any $1\leq i,j,k,\ell \leq n$. Denote by $T$ the torsion tensor of the Chern connection $\nabla^c$ of $g$. Its components under $e$ are given by $T(e_i, \overline{e}_j)=0$ and $T(e_i,e_k)=\sum_j T^j_{ik} e_j$. For our later proofs, we will need the following formula for Chern torsion components given by \cite[Lemma 1]{CZ}:

\begin{lemma} [\cite{CZ}]  \label{lemma1}
Given a Hermitian Lie algebra $({\mathfrak g}, J,g)$, let $e$ be a frame with dual coframe $\varphi$, and structure constants $C$ and $D$ be given by (\ref{CD1}). Then the Chern torsion components $T^j_{ik}$ under $e$ are given by
\begin{equation} \label{torsion}
 T^j_{ik}= - C^j_{ik} -  \sum_{\ell ,m=1}^n D^{m}_{\ell k} \,g^{\overline{m}j} g_{i \bar{\ell}} + \sum_{\ell ,m=1}^n  D^{m}_{\ell i} \,g^{\overline{m}j} g_{k \overline{\ell}} , \ \ \ \ \ \ \forall \ 1\leq i,j,k\leq n,
\end{equation}
where $g_{i\bar{j}}=\langle e_i, \overline{e}_j\rangle$ and $(g^{\overline{i}j})$ is the inverse matrix of $(g_{i\overline{j}})$. 
\end{lemma}

Now let us assume that ${\mathfrak g}$ is 2-step solvable, which means that ${\mathfrak g}$ itself is not abelian but the commutator ${\mathfrak g}'=[{\mathfrak g}, {\mathfrak g}]$ is abelian. From now on, we will ignore the trivial case and always assume that our ${\mathfrak g}$ is not abelian (note that for the abelian Lie algebra of even dimension, all Hermitian structures on it are K\"ahler and flat).  Let $(J,g)$ be a Hermitian structure on ${\mathfrak g}$. In \cite{FS, FS2}, Freibert and Swann systematically studied the Hermitian geometry of $2$-step solvable Lie algebras. Among other things, they give  characterizations to balanced and pluriclosed metrics on such Lie algebras, and confirmed Fino-Vezzoni Conjecture for all such Lie algebras that are of the pure types. Following their notations, we will write
\begin{equation} \label{UW}
{\mathfrak g}'_J = {\mathfrak g}' \cap J{\mathfrak g}', \ \ \     W=({\mathfrak g}' + J{\mathfrak g}')^{\perp}, \ \ \ \ U= ({\mathfrak g}'_J)^{\perp} \cap ({\mathfrak g}'+J{\mathfrak g}').
\end{equation}
Then we have the orthogonal decomposition of ${\mathfrak g}$ into $J$-invariant subspaces: ${\mathfrak g} = {\mathfrak g}'_J \oplus U\oplus W$. In the terminology of Freibert and Swann, the 2-step solvable Lie algebra ${\mathfrak g}$  is said to be of pure type I, II, or III, if the first, second, or third summand in the decomposition vanishes, namely when ${\mathfrak g}'\cap J{\mathfrak g}'=0$, ${\mathfrak g}'=J{\mathfrak g}'$, or ${\mathfrak g}'+J{\mathfrak g}'={\mathfrak g}$, respectively. Note that type II is exclusive to type I or III, but type I and III can overlap. Of course in general all three summands are present. Let us write
\begin{equation} \label{V}
V = ({\mathfrak g}'_J)^{\perp} \cap   {\mathfrak g}' ,  \ \ \ \  V' = V^{\perp} \cap U. 
\end{equation} 
Then we have the orthogonal decompositions ${\mathfrak g}'={\mathfrak g}'_J\oplus V$ and $U=V\oplus V'$, as well as the direct sum decomposition $U=V\oplus JV$ which may not be orthogonal. Following \cite{CZ}, let us introduce

\begin{definition}[{\bf admissible frames \cite{CZ}}] 
Let $({\mathfrak g}, J,g)$ be a $2$-step solvable Lie algebra equipped with a Hermitian structure. Write $W=({\mathfrak g}' + J {\mathfrak g}')^{\perp}$ and $V=({\mathfrak g}'_J)^{\perp} \cap {\mathfrak g}'$ as in (\ref{UW}) and (\ref{V}). Then a frame $e$ of ${\mathfrak g}$ is said to be {\em admissible} if $\{ e_1, \ldots , e_r\} $ is a unitary basis of $({\mathfrak g}'_J)^{1,0}$, and $\{ e_{s+1}, \ldots , e_n\}$ is a unitary basis of $W^{1,0}$, while $V$ is spanned by $ e_{\alpha} + \overline{e}_{\alpha}$ for $r\!+\!1\leq \alpha \leq s$. 
\end{definition}

Here we denoted by $2r$ the (real) dimension of ${\mathfrak g}'_J$ and by $2(n-s)$ the (real) dimension of $W$. Note that our admissible frames in general will not be unitary in the $U$ section, since $V$ may not be perpendicular to $JV$. As in \cite{CZ}, we will make the following convention on the range of indices: 
\begin{equation} \label{convention}
1\leq i,j,k,\ldots \leq r;  \  \ \ \ \ r\!+\!1\leq \alpha , \beta , \gamma , \ldots  \leq s;   \ \ \ \ \ s\!+\!1\leq a, b, c, \ldots  \leq  n.
\end{equation}

The following result was proved in \cite[Lemma 4]{CZ}, which gives the restriction on the structure constants for $2$-step solvable Lie algebras under any admissible frame:

\begin{lemma}[\cite{CZ}] \label{lemma2}
Let  $({\mathfrak g}, J,g)$ be a $2$-step solvable Lie algebra equipped with a Hermitian structure, and let $e$ be an admissible frame. Then 
\begin{equation} \label{res1}
 \left\{ \begin{split}  C^{\ast}_{ij}=C^{\alpha}_{\ast \ast} = C^{a}_{\ast \ast}= D^{\ast}_{a\ast} = D^i_{\ast j} = D^{\alpha}_{\ast j}=D^i_{\alpha \beta} = 0, \\
 C^j_{i\alpha}=-\overline{D^i_{j\alpha}}, \ \ \ \  C^{\ast}_{\alpha \beta} = \overline{ D^{\beta}_{\ast \alpha} }  - \overline{ D^{\alpha}_{\ast \beta} },  \hspace{2.0cm} \\  \overline{D^x_{\alpha y}} = - D^y_{\alpha x}, \ \ \ \ \forall \ 1\leq x,y\leq n.  \hspace{2.55cm}
 \end{split} \right.
\end{equation}
Here and from now on the indices $i,j,\alpha, \beta$ all follow the index range convention (\ref{convention}), while $\ast$ stands for an arbitrary integer between $1$ and $n$. 
\end{lemma}

For the sake of convenience, let us introduce the notation $C^Y_{XZ}=\sum_{i,j,k=1}^n X_iZ_k\overline{Y}_j C^j_{ik}$ for type $(1,0)$ vectors $X=\sum_{i=1}^n X_ie_i$, $Y=\sum_{i=1}^n Y_ie_i$ and $Z=\sum_{i=1}^n Z_ie_i$, where $e$ is a frame. Similarly, we will write that $D^Y_{XZ}=\sum_{i,j,k=1}^n X_iZ_k\overline{Y}_j D^j_{ik}$. Under such notations,  formula (\ref{CD2}) now takes the form
\begin{equation}
[X,Y] = \sum_{t=1}^n C_{XY}^t e_t , \ \ \ \ \ [X,\overline{Y} ] = \sum_{t=1}^n \big( \overline{D^X_{tY}} \,e_t -  D_{tX}^Y \overline{e}_t \big),   \label{CD3}
\end{equation}
for any type $(1,0)$ vectors $X$ and $Y$.

\vspace{0.3cm}

\section{Proof of Theorem \ref{thm3} }

In this section we will prove the main result of this article, Theorem \ref{thm3}, which states that any flat Hermitian Lie algebra is K\"ahler.

\begin{proof}[{\bf Proof of Theorem \ref{thm3}: }] Let $({\mathfrak g},J,g)$ be a Lie algebra with a Hermitian structure, so that the metric $g=\langle , \rangle$ is flat. By Milnor's Theorem \cite[Theorem 1.5]{Milnor}, ${\mathfrak g}$ is unimodular and is solvable of step at most 2. When ${\mathfrak g}$ is abelian, any Hermitian structure  $(J,g)$ on ${\mathfrak g}$ is K\"ahler flat. So it suffices to assume ${\mathfrak g}$ is non-abelian, hence it is 2-step solvable, and we have orthogonal direct sum ${\mathfrak g} = {\mathfrak h}\oplus {\mathfrak z}\oplus {\mathfrak g}'$ by the result of Barberis-Dotti-Fino, Theorem  \ref{thm2} stated in the introduction, where ${\mathfrak g}'$ is the commutator, ${\mathfrak z}$ is the center, and ${\mathfrak h}$ is an abelian subalgebra such that 
\begin{equation} \label{eq:skew} 
\langle [z, x], y\rangle + \langle [z,y], x\rangle =0, \ \ \ \ \forall \ z\in {\mathfrak h}, \ \ \forall \ x,y\in {\mathfrak g}. 
\end{equation}
Let $e$ be an admissible frame, so the structure constants $C$ and $D$ satisfy all the conditions in (\ref{res1}). For each $r+1\leq \alpha \leq s$, let us write $v_{\alpha} = \sum_{\beta =r+1}^s g^{\overline{\alpha} \beta} e_{\beta}$. Then the subspace $V'=V^{\perp}\cap U$ is given by
$$ V' = \mbox{span}_{\tiny \mathbb R}\{ \,\sqrt{-1} (v_{\alpha} - \overline{v}_{\alpha}) ; \, r< \alpha \leq s\}, $$
since obviously these $v_{\alpha}$ are linearly independent and $\langle v_{\alpha} - \overline{v}_{\alpha}, e_{\beta} + \overline{e}_{\beta} \rangle =0$ for any $\alpha$, $\beta$. By our notation (\ref{UW}) and (\ref{V}), we have ${\mathfrak g}'={\mathfrak g}'_J \oplus V$. Its orthogonal complement in ${\mathfrak g}$ is given by  $V'\oplus W = {\mathfrak h}\oplus {\mathfrak z}$, so in particular, 
$V'\oplus W$ is abelian and (\ref{eq:skew}) holds for any $z\in V'\oplus W$. The first part is equivalent to
\begin{equation} \label{eq:abelian}
[e_a, e_b]=[e_a, \overline{e}_b] = [v_{\alpha} - \overline{v}_{\alpha}, e_a]=  [v_{\alpha} - \overline{v}_{\alpha}, v_{\beta} - \overline{v}_{\beta}]=0 , \ \ \ \ \ \ \forall \ r< \alpha ,\beta \leq s, \ \forall\ s< a,b\leq n, 
\end{equation}
while the second part means
\begin{eqnarray} 
&& \langle [e_a,x],y\rangle + \langle [e_a,y],x\rangle =0, \ \ \ \ \ \ \ \forall \ s<a\leq n, \ \forall \ x,y \in {\mathfrak g}, \label{eq:eskew} \\
&& \langle [v_{\alpha} - \overline{v}_{\alpha} ,x],y\rangle + \langle [v_{\alpha} - \overline{v}_{\alpha},y],x\rangle =0, \ \ \ \ \ \ \forall \ r<\alpha \leq s, \ \forall \ x,y \in {\mathfrak g}. \label{eq:vskew}
\end{eqnarray}
Note that in (\ref{eq:eskew}) and (\ref{eq:vskew}) we may let $x$ and $y$ be any element in the complexification ${\mathfrak g}^{\mathbb C}= {\mathfrak g}^{1,0} \oplus \overline{{\mathfrak g}^{1,0}}  $. 
Using (\ref{CD3}), the formula (\ref{eq:abelian}) becomes
\begin{eqnarray} 
&& C_{ab}^{\ast} =D^a_{\ast b} = D^{\alpha }_{\ast a} =0, \ \ \ \ \ \ \ \ \forall\ r<\alpha \leq s, \ \forall \ s<a,b\leq n. \label{eq:abelian1} \\
&& C^{\, \ast}_{v_{\alpha}a} = - \overline{D^{\, a}_{\ast v_{\alpha}}} , \ \ \    C^{\, \ast}_{v_{\alpha}v_{\beta} } =  \overline{D^{v_{\alpha} }_{\ast v_{\beta}}}  - \overline{D^{v_{\beta} }_{\ast v_{\alpha}}} , \ \ \ \  \forall\ r<\alpha, \beta \leq s. \label{eq:abelian2}
\end{eqnarray}
For the last equality in (\ref{eq:abelian1}), we initially have $  D^{v_{\alpha }}_{\ast a} =0$ for any $\alpha $ and any $a$, but because the $(s-r)\times (s-r)$ matrix $(g_{\alpha \bar{\beta}})$ is non-degenerate, this is equivalent to  $D^{\alpha }_{\ast a} =0$ for any $\alpha $ and any $a$. Now by taking $(x, y)$ in (\ref{eq:eskew}) to be $(\overline{e}_x, \overline{e}_y)$ or $(e_x, \overline{e}_y)$, respectively, where $1\leq x,y\leq n$, we obtain
\begin{equation} 
 \sum_{t=1}^n \big( D^a_{tx}g_{y\bar{t}} + D^a_{ty}g_{x\bar{t}} \big) = \sum_{t=1}^n \big( C^t_{xa}g_{t\bar{y}} + D^y_{ta} g_{x\bar{t}} \big) =0, \ \ \ \ \ \ \  \forall \ s<a\leq n, \ \forall\ 1\leq x,y\leq n;  \label{eq:eskew1} 
\end{equation}
Similarly, (\ref{eq:vskew}) leads to
\begin{eqnarray} 
&&  \sum_{t=1}^n \big( D^{\alpha}_{tx}g_{y\bar{t}} + D^{\alpha}_{ty}g_{x\bar{t}} \big)   =0, \ \ \ \ \ \ \  \forall \ r<\alpha \leq s, \ \forall\ 1\leq x,y\leq n;  \label{eq:vskew1}  \\
&& \sum_{t=1}^n \{ (C^t_{\alpha x} + \overline{D^x_{t\alpha}}) g_{t\bar{y}} - ( \overline{ C^t_{\alpha y} } + D^y_{t\alpha}) g_{x\bar{t}} \} = 0, \ \ \ \ \ \ \  \forall \ r<\alpha \leq s, \ \forall\ 1\leq x,y\leq n;  \label{eq:vskew2} 
\end{eqnarray}
In particular, by taking $(x,y)=(i,j)$ or $(i,\alpha)$ in the first equality of (\ref{eq:eskew1}), we obtain
\begin{equation} 
 D^a_{ij}+D^a_{ji} \ = \ D^a_{i\alpha} + \sum_{\beta =r+1}^s \!D^a_{\beta i} g_{\alpha \bar{\beta}}  =0, \ \ \ \ \ \ \  \forall \ 1 \leq i,j\leq r, \ \forall \ r < \alpha \leq s, \ \forall  s<a\leq n.  \label{eq:eskew2} 
\end{equation}

\noindent {\em Claim 1:} $\ D^a_{\ast \ast}=0$, \ $\forall$ $s<a\leq n$. 

\begin{proof}[Proof of Claim 1:]  Let us take $i=j$ and $k=\ell =a$ in the third equation of (\ref{CCCD}). Since $D^a_{\ast b}=0$, we get
\begin{eqnarray*}
0 & = &  \sum_{t=1}^r  C^t_{ia}\overline{D^t_{ia}} -  \sum_{t=1}^s C^i_{ta}\overline{D^i_{ta}}   -  \sum_{t=1}^s D^a_{ti}\overline{D^a_{it}}  \\
& = &  \sum_{j=1}^r  \big( C^j_{ia}\overline{D^j_{ia}} -   C^i_{ja}\overline{D^i_{ja}} \big) - \sum_{\alpha =r+1}^s \big( C^i_{\alpha a}\overline{D^i_{\alpha a}} + D^a_{\alpha i}\overline{D^a_{i\alpha }} \big) - \sum_{j=1}^r D^a_{ji}\overline{D^a_{ij}}
\end{eqnarray*}
If we sum $i$ from $1$ to $r$, then the first two terms will cancel each other, so we just need to manage the rest. By taking $(x,y)=(\alpha ,i)$ in the second equation of (\ref{eq:eskew1}), we get $C^i_{\alpha a}=-\sum_{\beta} D^i_{\beta a}\,g_{\alpha \bar{\beta}}$. The second equality in (\ref{eq:eskew2}) gives us $D^a_{i\alpha } = - \sum_{\beta} D^a_{\beta i}\,g_{\alpha \bar{\beta}}$. Finally, the first equation of (\ref{eq:eskew2}) gives us $D^a_{ji}=-D^a_{ij}$. Plug all these into the above equality, and sum $i$ from $1$ to $r$, we get
$$  0 = 0 +   \sum_{i=1}^r\sum_{\alpha, \beta =r+1}^s D^i_{\beta a} \overline{ D^i_{\alpha a}} g_{\alpha \bar{\beta}} +  \sum_{i=1}^r\sum_{\alpha, \beta =r+1}^s D^a_{\alpha i}\overline{D^a_{\beta i}} g_{\beta \bar{\alpha}} + \sum_{i,j=1}^r |D^a_{ij}|^2 .
 $$
Since $(g_{\alpha \bar{\beta}})$ is positive definite, each term on the right hand side of the above equality is non-negative, hence we get $D^a_{ij}=D^a_{\alpha i}=D^a_{i\alpha}=0$. Since $D^a_{\ast b}=D^{\alpha}_{\ast a}=0$ and $D^x_{\alpha y} = - \overline{ D^y_{\alpha x}}$, together with (\ref{eq:eskew2}), we conclude that $D^a_{\ast \ast }=0$, and the claim is proved. 
\end{proof}

\noindent {\em Claim 2:} $\ D^{\alpha}_{\ast \ast}=0$, \ $\forall$ $r<\alpha \leq s$. 

\begin{proof}[Proof of Claim 2:] Since $D^{\alpha}_{\ast i} = D^{\alpha}_{\ast a} =0$, and by (\ref{eq:vskew1}) we have $D^{\alpha}_{i\beta}= -\sum_{\gamma} D^{\alpha}_{\gamma i} g_{\beta \bar{\gamma }} =0 $, so to prove Claim 2 it suffices to show $D^{\alpha}_{\beta \gamma }=0$ for all $\alpha$, $\beta$, $\gamma$ between $r+1$ and $s$.  Since $C^{\alpha}_{\ast \ast}=0$ and $C^{\ast}_{\beta \gamma } = \overline{D^{\gamma}_{\ast \beta} } - \overline{D^{\beta}_{\ast \gamma} } $, we have 
\begin{equation} \label{symmetry}
 D^{\gamma}_{\alpha \beta } = D^{\beta}_{\alpha \gamma } = -\overline{ D^{\gamma}_{\alpha \beta } }, \ \ \ \ \ \forall \ r<\alpha, \beta , \gamma \leq s,
 \end{equation}
where the last equality is due to the fact that $D^x_{\alpha y} = -\overline{ D^y_{\alpha x}} $ by the last line in (\ref{res1}),  in other words each $D^{\gamma}_{\alpha \beta } $ is pure imaginary. Also, (\ref{eq:vskew1}) says that
\begin{equation} \label{skew-symmetry}
 \sum_{\sigma =r+1}^s  \big( D^{\gamma}_{\sigma \beta } \,g_{\alpha \bar{\sigma}} + D^{\gamma}_{\sigma \alpha } \,g_{\beta \bar{\sigma}} \big) = 0, \ \ \ \ \ \forall \ r<\alpha, \beta , \gamma \leq s.
 \end{equation}
Now let $i,j,k,\ell$ be $\alpha, \beta , \gamma , \delta$, respectively in the last equation of (\ref{CCCD}), and using (\ref{symmetry}), we get
 \begin{equation} \label{commutativity}
 \sum_{\sigma =r+1}^s  \big( D^{\sigma}_{\gamma \alpha }  D^{\delta}_{\sigma \beta } - D^{\sigma}_{\gamma \beta }  D^{\delta}_{\sigma \alpha } \big) =0,  \ \ \ \ \ \forall \ r<\alpha, \beta , \gamma, \delta  \leq s.
 \end{equation}
Let us denote by $A_{\gamma}$ the real $(s-r)\times (s-r)$ matrix whose $(\alpha, \beta)$-th entry is $\sqrt{-1}D^{\beta}_{\alpha \gamma}$.  Again by (\ref{symmetry}), the equation (\ref{skew-symmetry}) now says that $GA_{\gamma}$ is skew-symmetric:
\begin{equation} \label{skew-symmetry2}
 GA_{\gamma} + \,^t\!(GA_{\gamma}) =0, \ \ \ \ \ \forall \ r<\gamma \leq s,
 \end{equation}
where $G$ is the positive definite $(s-r)\times (s-r)$ matrix $(g_{\alpha \bar{\beta}})$. Note that so far we only assumed that $e$ is an admissible frame for the 2-step solvable Lie algebra ${\mathfrak g}$, namely, $\{ e_i\}_{1\leq i\leq r}$ is a unitary frame for ${\mathfrak g}'_J$, $\{ e_a\}_{s<a\leq n}$ is a unitary frame for $W = ({\mathfrak g}'+J{\mathfrak g}')^{\perp}$, and $\{ e_{\alpha} \}_{r\leq \alpha \leq s}$ is any frame for the middle section $U=V\oplus JV$ such that $V=({\mathfrak g}'_J)^{\perp}\cap {\mathfrak g}'$ is spanned by $\{ e_{\alpha} + \overline{e}_{\alpha} \}_{r<\alpha \leq s}$. If we take any orthonormal basis $\{ \varepsilon_{\alpha}\}_{r<\alpha \leq s}$ for $V$, and take $e_{\alpha} = \frac{1}{\sqrt{2}} ( \varepsilon_{\alpha} - \sqrt{-1} J\varepsilon_{\alpha})$ for each $\alpha$, then for this particular admissible frame $e$ we have
$$ \langle e_{\alpha}, \overline{e}_{\beta} \rangle = \frac{1}{2} \langle  \varepsilon_{\alpha} - \sqrt{-1} J\varepsilon_{\alpha},  \varepsilon_{\beta} + \sqrt{-1} J\varepsilon_{\beta} \rangle = \delta_{\alpha \beta} + \sqrt{-1} H_{\alpha \beta}, $$
where $H_{\alpha \beta } = \langle \varepsilon_{\alpha} , J\varepsilon_{\beta}\rangle = - \langle \varepsilon_{\beta} , J\varepsilon_{\alpha}\rangle$ is real. So for our frame $e$ the matrix of metric $G=I+\sqrt{-1}H$ where $H$ is real and skew-symmetric. Now if we take the real part in equation (\ref{skew-symmetry2}), we see that each $A_{\gamma}$ is skew-symmetric, namely, 
\begin{equation} \label{skew-symmetry3}
 D^{\beta}_{\alpha \gamma} = - D^{\alpha}_{\beta \gamma}, \ \ \ \ \ \forall \ r<\alpha , \beta , \gamma \leq s.
 \end{equation}
By (\ref{symmetry}), this could also be written as $ D^{\gamma}_{\alpha \beta} = - D^{\gamma}_{\beta \alpha}$. Therefore, the commutativity equation (\ref{commutativity})  can be written as 
\begin{equation} \label{commutativity2}
 \sum_{\sigma =r+1}^s  \big( D^{\sigma}_{\alpha \gamma }  D^{\delta}_{\beta\sigma  } - D^{\sigma}_{\beta \gamma  }  D^{\delta}_{\alpha\sigma  } \big) =0,  \ \ \ \ \ \forall \ r<\alpha, \beta , \gamma, \delta  \leq s,
 \end{equation}
 that is, $[B_{\alpha} , B_{\beta}]=0$, where $B_{\alpha} = (\sqrt{-1}D^{\gamma}_{\alpha \beta})$ is the real, symmetric $(s-r)\times (s-r)$ matrix. Denote by ${\mathcal B}$ the linear space spanned by $B_{\alpha}$ for all $r<\alpha \leq s$. The matrices in ${\mathcal B}$ can be simultaneously diagonalized by orthogonal matrices, in other words, by an orthogonal change of the basis $\{ \varepsilon_{\alpha}\}_{r<\alpha \leq s}$, we may assume that all $B_{\alpha}$ are diagonal, thus
 $$ D^{\gamma}_{\alpha \beta} = \sqrt{-1} \lambda_{\alpha \beta} \delta_{\beta \gamma}, \ \ \ \ \forall\ r<\alpha, \beta, \gamma \leq s, $$
for some real constants $\lambda_{\alpha \beta} $. So $ D^{\gamma}_{\alpha \beta}=0$ when $\beta \neq \gamma$, while for $D^{\beta}_{\alpha \beta}$, by (\ref{skew-symmetry3}) it equals to $-D^{\alpha}_{\beta \beta}$, which is zero if $\alpha \neq \beta$, and it is also zero when $\alpha=\beta$ as in this case we have $D^{\beta}_{\beta \beta} = -  D^{\beta}_{\beta\beta}$. We have therefore proved that all $D^{\gamma}_{\alpha \beta} =0$, hence the claim.
\end{proof}

By the above two claims,  Lemma \ref{lemma2}, and the formula (\ref{eq:abelian1}), (\ref{eq:abelian2}), we know that the only possibly non-trivial components for $C$ and $D$ are $C^j_{i\alpha}$, $C^j_{ia}$, $D^j_{i\alpha}$, and $D^j_{ia}$. Also, the second equation in (\ref{eq:eskew1}) gives us $C^j_{ia}=-D^j_{ia}$, while a combination of (\ref{eq:vskew2}) and the fact $C^j_{ia}= -\overline{D^i_{j\alpha}}$ from (\ref{res1}) gives us $C^j_{i\alpha }=-D^j_{i\alpha }$. 

Now we are finally ready to show that the metric $g$ must be K\"ahler, which is equivalent to the condition that the Chern torsion tensor vanishes. To check this, we will use the formula (\ref{torsion}) given in Lemma \ref{lemma1}. First of all, since $C^{\alpha}_{\ast\ast}=D^{\alpha}_{\ast\ast }=0$, we get $T^{\alpha}_{\ast \ast}=0$ for any $r+1\leq \alpha \leq s$. Similarly,  since $C^{a}_{\ast\ast}=D^{a}_{\ast\ast }=0$, we get $T^{a}_{\ast \ast}=0$ for any $s+1\leq a\leq n$. It remains to check that $T^i_{xy }=0$ for any $x$ and $y$. We have
$$ C^{\ast}_{ij} = C^{\ast}_{\alpha \beta} = C^{\ast}_{ab} =C^{\ast}_{\alpha a}=0, \ \ \ D^{\ast}_{ij} = D^{\ast}_{\alpha \beta} = D^{\ast}_{ab} =D^{\ast}_{\alpha a}=0. $$
Since $T^i_{xy}=-T^i_{yx}$, we only need to check for $T^i_{j\alpha}$ and $T^i_{ja}$. We have
$$ T^i_{j\alpha} = - C^i_{j\alpha} - D^i_{j\alpha} + \sum_{\beta =r+1}^s \! D^i_{\beta j} \, g_{\alpha \bar{\beta}} =  - (C^i_{j\alpha} + D^i_{j\alpha}) = 0.$$
Similarly, 
$$ T^i_{ja} = - C^i_{ja} - D^i_{ja} +  D^i_{a j}  =  - (C^i_{ja} + D^i_{ja}) = 0.$$  
Thus the Chern torsion vanishes, and the metric $g$ is K\"ahler. This has completed the proof of Theorem \ref{thm3}.
\end{proof}

In Theorem \ref{thm3}, the flatness of $g$ means Riemannian flat, that is, the Levi-Civita (Riemannian) connection has vanishing curvature. It becomes an interesting question if one replaces the Levi-Civita connection by the Chern connection or the Bismut connection. In other words, one may ask the following

\begin{question}
Besides complex Lie groups, what kind of unimodular Lie groups with left-invariant Hermitian structure will be Chern flat? Besides Samelson spaces, what unimodular Lie groups with left-invariant Hermitian structure will be Bismut flat?
\end{question}

On a complex Lie group, any left-invariant metric compatible with the complex structure is always Chern flat, namely, with all $D^j_{ik}=0$. Conversely, the classic Boothby's Theorem \cite{Boothby} states that the universal cover of any compact Chern flat manifold $(M^n,g)$  is biholomorphic to a complex Lie group with a left-invariant metric. However, it may happen (and there are indeed such examples) that the universal cover of $(M^n,g)$ is also biholomorphic to a Lie group $G$ with a left-invariant Hermitian structure, which is Chern flat yet the group $G$ is not a complex Lie group. So the above question meant to classify all such Hermitian Lie groups. 

Similarly, a Samelson space means a connected, simply-connected Lie group $G$ with a bi-invariant metric and a compatible left-invariant complex structure. They are known to be Bismut flat, and Milnor's Lemma \cite{Milnor} says that such Lie groups are exactly the product of the vector group ${\mathbb R}^k$ with  compact semisimple Lie groups. In \cite{WYZ}, Wang, Yang and the second named author proved that the universal cover of any complete Bismut flat manifold is biholomorphic to a Samelson space. Again since different Lie groups with Hermitian structure could be biholomorphic to each other yet non-isomorphic as groups, the above question intends to classify all Bismut flat Hermitian Lie algebras other than the Samelson ones. 

While it might be difficult to answer the above questions in their full generality, it might be a more feasible task if we restrict ourselves to only unimodular 2-step solvable groups. Hopefully similar techniques used in this article could lead to an answer.

\vspace{0.5cm}

\noindent\textbf{Acknowledgments.} {We would like to thank Haojie Chen, Shuwen Chen, Lei Ni, Xiaolan Nie, Kai Tang, Bo Yang, Yashan Zhang, and Quanting Zhao for their interest and helpful discussions.}


\end{document}